\documentclass{amsart}
\usepackage{amsmath,amsfonts,amsthm,enumerate,color}
\newcommand{\End}{\text{\rm End}}
\newcommand{\thmref}[1]{Theorem~\ref{#1}}
\newcommand{\propref}[1]{Proposition~\ref{#1}}

\newcommand{\lemref}[1]{Lemma~\ref{#1}}
\newcommand{\eqnref}[1]{~(\ref{#1})}

\newcommand{\germ}{\mathfrak}
\newtheorem{thm}{Theorem}[section]
\newtheorem{lem}[thm]{Lemma}

\theoremstyle{definition}

\newtheorem{prop}[thm]{Proposition}

\subjclass{Primary 17B67, 81R10}

\theoremstyle{rem}
\numberwithin{equation}{section}

\usepackage{amssymb}
\begin{document}
\title{Free Field Realizations of the Elliptic Affine Lie Algebra $\mathfrak{sl}(2, R)  \oplus\left( \Omega_R/dR\right)$}
\author{Andr\'e Bueno}\author{ Ben Cox}\author{Vyacheslav Futorny}
\begin{abstract}  In this paper we construct two free field
 realizations of  the elliptic affine Lie algebra
$\mathfrak{sl}(2,R)\oplus \Omega_R/dR$ where $R=\mathbb
C[t,t^{-1},u\,|\, u^2=t^3 - 2b t^2 + t]$.  The first realization
provides an analogue of Wakimoto's construction for Affine Kac-Moody algebras, but in the setting of the elliptic affine Lie algebra.
The second realization gives new type of representations analogous to
Imaginary Verma modules in the Affine setting.
\end{abstract}
\keywords{Wakimoto Modules, Elliptic Affine Algebras, Affine Lie
Algebras, Fock Spaces}
\address{Department of Mathematics \\
University of Charleston \\
66 George Street  \\
Charleston SC 29424, USA\\ (second author)}\email{coxbl@cofc.edu}
\urladdr{http://math.cofc.edu/faculty/cox/papers/}
\address{Institute of Mathematics, University of S\~ao Paulo\\
S\~ao Paulo, 05315-970\\
 Brazil (first and third authors)}
\email{futorny@ime.usp.br}
\maketitle
\section{Introduction}
Elliptic affine Lie algebras are a particular family of
Krichever-Novikov Lie algebras related to Riemann surfaces
\cite{KN1}, \cite{KN2}. These algebras were introduced by
Krichever and Novikov in their study of string theory in Minkowski
space.

The theory of highest weight modules for these algebras was developed
by O. Sheinman in \cite{MR1082031}, \cite{MR1328538} and \cite{MR1328538}, see
also \cite{MR1027505}. Since Krichever-Novikov algebras are
quasi-graded, their representation theory is quite different from
the standard representation theory of Kac-Moody algebras. For
instance, irreducible highest weight modules can have a $2$-dimensional
subspace of highest weight vectors.

The goal of present paper is to obtain free field realizations
for the elliptic affine Lie algebra associated with
$\mathfrak{sl}(2, \mathbb C)$.

It is known from the work of Kassel and Loday (see \cite{MR694130}, and \cite{MR772062}) that if $R$ is a commutative algebra and $\mathfrak g $ is a simple Lie algebra, both defined over the complex numbers, then the universal central extension $\hat{{\mathfrak g}}$ of
$L({\mathfrak g})={\mathfrak g}\otimes R$ is the vector space $L({\mathfrak g})\oplus \Omega_R^1/dR$ where $\Omega_R^1/dR$ is the space of K\"ahler differentials modulo exact forms (see \cite{MR772062}).  The vector space $\hat{{\mathfrak g}}$ is made into a Lie algebra by defining
$$
[x\otimes f,y\otimes g]:=[xy]\otimes fg+(x,y)\overline{fdg},\quad [x\otimes f,\omega]=0
$$
for $x,y\in\mathfrak g$, $f,g\in R$,  here $\omega\in \Omega_R^1/dR$ and $(-,-)$ denotes the Killing form on $\mathfrak g$ and $\overline{a}$ denotes the image of $a\in\Omega^1_R$ in the quotient $\Omega^1_R/dR$.     A somewhat vague (due to the imprecision in the choice of $R$ and hence imprecision in the description of the basis of $\Omega_R/dR$) but natural question is whether there exist free field or Wakimoto type realizations of these algebras.    The answer is well known from the work of M. Wakimoto if $R$ is the ring of Laurent polynomials in one variable (see \cite{W}). We  solve this problem in the setting where $\mathfrak g=\mathfrak{sl}(2,\mathbb C)$ and $R=\mathbb C[t,t^{-1},u|u^2=t^3 - 2b t^2 + t]$ is the elliptic algebra.
Related work on realizations of the universal central extension of $\mathfrak{sl}(2,\mathbb C)\otimes R$ can be found in the following papers; \cite{MR1303073}, \cite{MR2373448}, \cite{c2}, \cite{MR2065632}, \cite{MR89k:17016}, \cite{MR92f:17026}, \cite{MR92d:17025}, \cite{JK}, \cite{MR1082031}, \cite{W}.

Before we begin we would like to mention a little genesis of elliptic affine algebras.   In Kazhdan and Lusztig's explicit study of the
tensor structure of modules for affine Lie algebras (see
\cite{MR1186962} and \cite{MR1104840}) the ring of functions
regular everywhere except at a finite number of points appears
naturally.    This algebra M.\,Bremner gave the name {\it $n$-point
algebra}.   On the other hand an elliptic algebra is an algebra of functions on an elliptic curve of genus one which may have poles at two points.  This article deals with the particular example of a fixed nonsingular compact complex algebraic curve of genus 1 which Bremner denotes by $\Sigma$ in \cite{MR1303073}. This curve $\Sigma$ can be represented as a quotient of the complex plane $\mathbb C$ by the lattice $\Lambda=\mathbb Z\oplus \mathbb Z\lambda$ with positive imaginary part of $\lambda$.
Let $R$ denote ring of meromorphic functions on $\Sigma$ which are holomorphic outside of the set
$\displaystyle{\{0,\frac{1}{2}(1+\lambda)\}}$.
Then Bremner has shown that $R\cong S_b:=\mathbb C[t,t^{-1},u|u^2=t^3 - 2b t^2 + t]$ where $b\in \mathbb C$ is given below. As the latter,
being $\mathbb Z_2$-graded, is  more immediately amendable to the
theatrics of conformal field theory,
 we choose to work with $S_b$ instead.  Moreover, Bremner has given an
explicit  description of the universal central extension of $\mathfrak
g\otimes R$, in terms of Pollaczek polynomials.    His description is recapitulated in what follows.

Our main result, \thmref{mainthm}, provides a natural free field
realization in terms of a $\beta$-$\gamma$-system and the
oscillator algebra of the elliptic affine Lie algebra when
$\mathfrak g=\mathfrak{sl}(2,\mathbb C)$.   Just as in the case of
intermediate Wakimoto modules defined in \cite{MR2271362}, there two
different realizations depending on two different normal
orderings. The first realization is analogous to the construction of
Wakimoto modules for Affine Lie algebras while the second one
provides new  modules our algebra analogous to Imaginary Verma
modules \cite{MR95a:17030}.

\dedicatory{\it This paper is dedicated to Stephen Berman.}

\section{The ring of elliptic functions $\mathbb C[t,t^{-1},u|u^2=t^3 - 2b t^2 + t]$.}

Fix a nonsingular compact complex algebraic curve of genus 1 which we denote by $\Sigma$. This curve $\Sigma$ can be represented as a quotient of the complex plane $\mathbb C$ by the lattice $\Lambda=\mathbb Z\oplus \mathbb Z\lambda$ with positive imaginary part of $\lambda$.
Let $R$ denote ring of meromorphic functions on $\Sigma$ which are holomorphic outside of the set
$\displaystyle{\{0,\frac{1}{2}(1+\lambda)\}}$.

Set $\displaystyle{m=\wp\left(\frac{1}{2}(1+\lambda)\right)}$ where $\wp$ is the Weierstrass $\wp$ function:
$$
\wp(z):=z^{-2}+\sum_{0\neq \omega\in\Lambda}\left(\frac{1}{(z-\omega)^2}-\frac{1}{\omega^2}\right)
$$
Then
$$
\wp'(z)^2=4\wp(z)^3-g_2\wp(z)-g_3,\quad g_2=60\sum_{0\neq \xi \in \Lambda}\xi^{-4},\quad
g_3=140\sum_{0\neq \xi \in \Lambda}\xi^{-6}.
$$

\begin{prop}[\cite{MR1303073}, Prop. 4.1.]
We have
$$
R\cong  \mathbb C[t,t^{-1},u|u^2=t^3 - 2b t^2 + t]
$$
where $b$ is some constant determined by $m$.
\end{prop}

\section{The Universal Central Extension of $\mathfrak g\otimes R$}
We recall the realization of the universal central extension of $L({\mathfrak g})=\mathfrak g\otimes R$.
\begin{thm}[\cite{MR1303073}, cf. Theorem 3.4]  The space $\Omega_R^1/dR$ has a basis
$$
\omega_0:=\overline{t^{-1}dt},\quad \omega_-:=\overline{t^{-2}u\,dt},\quad \omega_+:=\overline{t^{-1}u\,dt}.
$$
\end{thm}
There exists an automorphism $\tau$ of $R$ given by
\begin{equation}
  \tau(t)=t^{-1},\quad \tau(u)=t^{-2}u
\end{equation}
which induces an automorphism on $\Omega_R^1/dR$.   This induced automorphism is simply the negative of the identity map.

We will now give Bremner's Fourier mode description of the relations satisfied by the basis elements $x\otimes t^n$, $x\otimes t^nu$, $\omega_0$, $\omega_\pm$ of $\hat{{\mathfrak g}}$.  First  recall the {\it Pollaczek polynomials} $P_k(b)=P_k^\lambda(b;\alpha,\beta,\gamma)$,
$\alpha,\beta,\gamma\in \mathbb C$ (see \cite{MR0035868}), which are defined by the recursion
\begin{equation}
(k+\gamma)P_k(b)=2[(k+\lambda+\alpha+\gamma-1)b+\beta]P_{k-1}(b)-(k+2\lambda +\gamma-2)P_{k-2}(b).
\end{equation}
\begin{lem}[\cite{MR1303073}, Lemma 4.4]
Consider the sequence of polynomials $p_k(b)$, and $q_k(b)$, defined by
$$
\overline{t^{k-2}u\, dt}=p_k(b)\overline{t^{-1}u\,dt}+q_k(b)\overline{t^{-2}u\,dt}.
$$
These polynomials are Pollaczek polynomials for the parameters $\lambda=-1/2$, $\alpha=0$, $\beta=-1$, $\gamma=1/2$ together with the initial conditions
$$
p_0(b)=0,\quad p_1(b)=1,\quad q_0(b)=1,\quad q_1(b)=0.
$$
\end{lem}
If we set
$$
Q(x,b):=\sum_{n=0}^\infty q_n(b)x^n,\quad P(x,b):=\sum_{n=0}^\infty p_n(b)x^n,
$$
then it is straightforward to show that
$$
x(x^2-2bx-1)\frac{dQ}{dx}+\left[(2\lambda+\gamma)x^2-2x((\lambda+\alpha+\gamma)b+\beta)+\gamma\right]Q =\gamma,
$$
and
$$
x(x^2-2bx-1)\frac{dP}{dx}+\left[(2\lambda+\gamma)x^2-2x((\lambda+\alpha+\gamma)b+\beta)+\gamma\right] P=(1+\gamma)x.
$$

Set $\alpha_{\pm}:=b\pm \sqrt{b^2-1}$, $A_{\pm}:=\lambda \mp (\alpha b+\beta)/\sqrt{b^2-1}$.  Solving these differential equations above we get
$$
Q(x,b)=\int_0^x\frac{\gamma \xi^{\gamma-1}(\xi-\alpha_1)^{A_+ -1}(\xi-\alpha_2)^{A_- -1}}{x^\gamma (x-\alpha_1)^{A_+}(x-\alpha_2)^{A_-}}\,d\xi
$$
and
$$
P(x,b)=\int_0^x\frac{(1+\gamma)\xi^{\gamma}(\xi-\alpha_1)^{A_+ -1}(\xi-\alpha_2)^{A_- -1}}{x^\gamma(x-\alpha_1)^{A_+}(x-\alpha_2)^{A_-}}\,d\xi.
$$

\begin{thm}[\cite{MR1303073}, Theorem 4.6] \label{Bremnersthm}
The elliptic Lie algebra $\hat{{\mathfrak g}}$ has a $\mathbb Z/2\mathbb Z$-grading where
\[
\hat{{\mathfrak g}}^0=\mathfrak g\otimes \mathbb C[t,t^{-1}]\oplus \mathbb C\omega_0,\quad
\hat{{\mathfrak g}}^1=\mathfrak g\otimes \mathbb C[t,t^{-1}]u\oplus \mathbb C\omega_-\oplus \mathbb C\omega_+.
\]
For $x,y\in\mathfrak g$, the commutation relations defining $\hat{{\mathfrak g}}$ are
\begin{align*}
[x\otimes t^i,y\otimes t^j]&=[xy]\otimes t^{i+j}+\delta_{i+j,0}j\omega_0  \\
[x\otimes t^{i-1}u,y\otimes t^{j-1}u]&=[xy]\otimes\left(t^{i+j-1}-2bt^{i+j}+t^{i+j+1}\right) \\
 &\hskip 20pt +
(x,y)\omega_0\left(-2j b\delta_{i+j,0}+\frac{1}{2}(j-i)(\delta_{i+j,-1}+\delta_{i+j,1})\right)  \\
[x\otimes t^{i-1}u,y\otimes t^{j}]&=[xy]\otimes t^{i+j-1} u+(x,y)j\Big(p_{|i+j| }(b) \omega_+ +q_{|i+j|}(b)\omega_-\Big),
\end{align*}
for $i,j\in\mathbb Z$.
 In addition the elements $\omega_0,\omega_\pm $ are central.
\end{thm}


\subsection{Formal Distributions}
We introduce some notation that will simplify later arguments.
This notation follows roughly \cite{MR99f:17033} and \cite
{MR2000k:17036}:  The {\it formal delta function}
$\delta(z/w)$ is the formal distribution
$$
\delta(z/w)=z^{-1}\sum_{n\in\mathbb Z}\left(\frac{z}{w}\right)^n
$$
Given a vector space $V$, for any sequence of elements $\{a_{m}\}_{m\in
\mathbb Z}$ in the ring $\End (V)$,  the
formal distribution
\begin{align*}
a(z):&=\sum_{m\in\mathbb Z}a_{m}z^{-m-1}
\end{align*}
is called a {\it field}, if for any $v\in V$, $a_{m}v=0$ for $m\gg0$.
If $a(z)$ is a field, then we set
\begin{align}\label{usualnormalordering}
    a(z)_-:&=\sum_{m\geq 0}a_{m}z^{-m-1},\quad\text{and}\quad
   a(z)_+:=\sum_{m<0}a_{m}z^{-m-1}.
\end{align}
 The {\it normal ordered product} of two distributions
$a(z)$ and
$b(w)$ (and their coefficients) is
defined by
\begin{equation}\label{normalorder}
\sum_{m\in\mathbb Z}\sum_{n\in\mathbb
Z}:a_mb_n:z^{-m-1}w^{-n-1}=:a(z)b(w):=a(z)_+b(w)+b(w)a(z)_-.
\end{equation}

Note that while $:a^1(z_1)\cdots a^m(z_m):$
is always defined as a formal series, we will only define $:a(z)
b(z):=\lim_{w\to z}:a(z)b(w):$
for certain pairs
$(a(z),b(w))$.

Then one defines recursively
\[
:a^1(z_1)\cdots a^k(z_k):=:a^1(z_1)\left(:a^2(z_2)\left(:\cdots
:a^{k-1}(z_{k-1}) a^k(z_k):\right)\cdots
:\right):,
\]
while normal ordered product
\[
:a^1(z)\cdots
a^k(z):=\lim_{z_1,z_2,\cdots, z_k\to
z} :a^1(z_1)\left(:a^2(z_2)\left(:\cdots :a^{k-1}(z_{k-1})
a^k(z_k):\right)\cdots
\right):
\]
will only be defined for certain $k$-tuples $(a^1,\dots,a^k)$.

Let
\begin{equation}\label{contraction}
\lfloor
a(z)b(w)\rfloor=\lfloor
ab\rfloor=a(z)b(w)-:a(z)b(w):= [a(z)_-,b(w)],
\end{equation}
where
$\lfloor
a(z)b(w)\rfloor$ denotes the {\it contraction} of any two formal distributions
$a(z)$ and $b(w)$.
\begin{thm}[Wick's Theorem, see \cite{MR85g:81096}, \cite{MR99m:81001} or
\cite{MR99f:17033} ]  Let  $a^i(z)$ and $b^j(z)$ be formal
distributions with coefficients in the associative algebra
 $\End(\mathbb C[\mathbf x]\otimes \mathbb C[\mathbf y])$,
 satisfying
\begin{enumerate}
\item $[ \lfloor a^i(z)b^j(w)\rfloor ,c^k(x)_\pm]=[ \lfloor
a^ib^j\rfloor ,c^k(x)_\pm]=0$, for all $i,j,k$ and
$c^k(x)=a^k(z)$ or
$c^k(x)=b^k(w)$.
\item $[a^i(z)_\pm,b^j(w)_\pm]=0$ for all $i$ and $j$.
\item The products
$$
\lfloor a^{i_1}b^{j_1}\rfloor\cdots
\lfloor a^{i_s}b^{i_s}\rfloor:a^1(z)\cdots a^M(z)b^1(w)\cdots
b^N(w):_{(i_1,\dots,i_s;j_1,\dots,j_s)}
$$
have coefficients in
$\End(\mathbb C[\mathbf x]\otimes \mathbb C[\mathbf y])$ for all subsets
$\{i_1,\dots,i_s\}\subset \{1,\dots, M\}$, $\{j_1,\dots,j_s\}\subset
\{1,\cdots N\}$. Here the subscript
${(i_1,\dots,i_s;j_1,\dots,j_s)}$ means that those factors $a^i(z)$,
$b^j(w)$ with indices
$i\in\{i_1,\dots,i_s\}$, $j\in\{j_1,\dots,j_s\}$ are to be omitted from
the product
$:a^1\cdots a^Mb^1\cdots b^N:$ and when $s=0$ we do not omit
any factors.
\end{enumerate}
Then
\begin{align*}
:&a^1(z)\cdots a^M(z)::b^1(w)\cdots
b^N(w):= \\
  &\sum_{s=0}^{\min(M,N)}\sum_{i_1<\cdots<i_s,\,
j_1\neq \cdots \neq j_s}\lfloor a^{i_1}b^{j_1}\rfloor\cdots
\lfloor a^{i_s}b^{i_s}\rfloor
:a^1(z)\cdots a^M(z)b^1(w)\cdots
b^N(w):_{(i_1,\dots,i_s;j_1,\dots,j_s)}.
\end{align*}
\end{thm}

For $m=i-\frac{1}{2}$, $i\in\mathbb Z+\frac{1}{2}$ and $x\in\mathfrak g$, define $x_{m+\frac{1}{2}}:=x\otimes t^{i-\frac{1}{2}}u$ and $x_m:=x\otimes t^m$.  We set
\begin{align*}
x^1(z)
&:=\sum_{m\in\mathbb Z}x_{m+\frac{1}{2}}z^{-m-1},\quad x(z):=\sum_{m\in\mathbb Z}x_{m}z^{-m-1}.
\end{align*}
Then the relations in \thmref{Bremnersthm} can be rewritten as
\begin{align*}
[x(z),y(w)]
&= [xy](w)\delta(z/w)-(x,y)\omega_0\partial_w\delta(z/w) \\ \\
[x^1(z),y^1(w)]
&= w(w^2 -2bw+1)\left([x,y](w)\delta(z/w) -(x,y)\omega_0\partial_w\delta(z/w)\right)-(x,y)\frac{\omega_0}{2}\partial_w(w(w^2 -2bw+1))  \delta(z/w)  \\ \\
[x^1(z),y(w)]
&=[x,y]^1(w)\delta(z/w) \\
&\quad-(x,y) \left((P(w^{-1},b)+P(w,b))\omega_++(Q(w^{-1},b)+Q(w,b)-2)\omega_-\right) \delta(z/w)\\
&\quad-(x,y) w \partial_w\left(\left((P(w^{-1},b)+P(w,b))\omega_++(Q(w^{-1},b)+Q(w,b)-2)\omega_-\right)\delta(z/w)\right)\\  \\
[x(z),y^1(w)]
&=[x,y]^1(w)\delta(w/z)  \\
&\quad -(x,y)\left((P(w^{-1},b)+P(w,b))\omega_++(Q(w^{-1},b)+Q(w,b)-2)\omega_-\right)w \partial_w\delta(z/w) \\
 \end{align*}

 \section{ Oscillator algebras}
 \subsection{The $\beta-\gamma$ system}
 Let $\hat{\mathfrak a}$ be the infinite dimensional oscillator algebra with generators $a_n,a_n^*,a^1_n,a^{1*}_n,\,n\in\mathbb Z$ together with $\mathbf 1$ satisfying the relations
\begin{gather*}
[a_n,a_m]=[a_m,a_n^1]=[a_m,a_n^{1*}]=[a^*_n,a^*_m]=[a^*_n,a^1_m]=[a^*_n,a^{1*}_{m}]=0,\\
[a_n^{1},a_m^{1}]=[a_n^{1*},a_m^{1*}]=0=[\mathfrak a,\mathbf 1], \\
[a_n,a_m^*]=\delta_{m+n,0}\mathbf 1=[a^1_n,a_m^{1*}].
\end{gather*}
For  $c=a,a^1$ and respectively $X=x,x^1$ with $r=0$ or $r=1$, we define
a representation
$\rho_r:\hat{\mathfrak a}\to \mathfrak{gl}(\mathbb C[x_n,x_n^1\,|\,n\in\mathbb Z])$
of $\hat{\mathfrak a}$ in the Fock space $\mathbb C[x_n,x_n^1\,|\,n\in\mathbb Z])$
by
\begin{align*}
\rho_r( c_{m}):&=\begin{cases}
  \partial/\partial
X_{m}&\quad \text{if}\quad m\geq 0,\enspace\text{and}\enspace  r=0
\\ X_{m} &\quad \text{otherwise},
\end{cases}
 \\
\rho_r(c_{m}^*):&=
\begin{cases}X_{-m} &\enspace \text{if}\quad m\leq
0,\enspace\text{and}\enspace r=0\\ -\partial/\partial
X_{-m}&\enspace \text{otherwise}. \end{cases}
\end{align*}
and $\rho_r(\mathbf 1)=1$.
These two representations can be constructed using induction:
For $r=0$ the representation
$\rho_0$ is the
$\hat{\mathfrak a}$-module generated by $1=:|0\rangle$, where
$$
a_{m}|0\rangle=a^1_{m}|0\rangle=0,\quad m\geq  0,
\quad a_{m}^*|0\rangle= a_{m}^{1*}|0\rangle=0,\quad m>0.
$$
For $r=1$ the representation
$\rho_1$ is the
$\hat{\mathfrak a}$-module generated by $1=:|0\rangle$, where
$$
\quad a_{m}^*|0\rangle= a_{m}^{1*}|0\rangle=0,\quad m\in\mathbb Z.
$$
If we write
\begin{align*}
 \alpha(z):&=\sum_{n\in\mathbb Z}a_nz^{-n-1},\quad  &\alpha^*(z):=\sum_{n\in\mathbb Z}a_n^*z^{-n}  \\
 \alpha^1(z):&=\sum_{n\in\mathbb Z}a_n^1z^{-n-1},\quad  &\alpha^{1*}(z):=\sum_{n\in\mathbb Z}a_n^{1*}z^{-n}
 \end{align*}
then
\begin{align*}
[\alpha(z),\alpha(w)]&=[\alpha^*(z),\alpha^*(w)]=[\alpha^{1}(z),\alpha^{1}(w)]=[\alpha^{1*}(z),\alpha^{1*}(w)]=0  \\
[\alpha(z),\alpha^*(w)]&=[\alpha^1(z),\alpha^{1*}(w)]
    =\mathbf 1\delta(z/w).
\end{align*}
Observe that $\rho_1(a(z))$ and
$\rho_1(a^1(z))$ are not fields whereas $\rho_r(a^*(z))$ $\rho_r(a^{1*}(z))$
are always the fields.
Corresponding to these two representations there are two possible normal orderings:  For $r=0$ we use the usual normal ordering given by \eqnref{usualnormalordering} and for $r=1$ we define the {\it natural normal ordering} to be
\begin{alignat*}{2}
\alpha(z)_+&=\alpha(z),\quad &\alpha(z)_-&=0 \\
\alpha^1(z)_+&=\alpha^1(z),\quad &\alpha^1(z)_-&=0 \\
\alpha^*(z)_+&=0,\quad &\alpha^*(z)_-&=\alpha^*(z), \\
\alpha^{1*}(z)_+&=0,\quad &\alpha^{1*}(z)_-&=\alpha^{1*}(z) ,
\end{alignat*}
This means in particular that for $r=0$ we get
\begin{align}
\lfloor \alpha\alpha^*\rfloor
&
=\delta_-(z/w)
=\
\,\iota_{z,w}\left(\frac{1}{z-w}\right)\\
\lfloor \alpha^*\alpha\rfloor
&
=-\delta_+(w/z)=\,\iota_{z,w}\left(\frac{1}{w-z},
\right)
\end{align}
and for $r=1$
\begin{align}
\lfloor \alpha\alpha^*\rfloor
&=[\alpha(z)_-,\alpha^*(w)]=0 \\
\lfloor \alpha^*\alpha\rfloor
&
=[\alpha^*(z)_-,\alpha(w)]
=- \delta(w/z),
\end{align}
where similar results hold for $\alpha^1$.
Notice that in both cases we have
$$
\lfloor \alpha(z)\alpha^*(w)\rfloor-\lfloor\alpha^*(w) \alpha(z)\rfloor=\delta(z/w)
$$
We will also need the following two
results.
\begin{thm}[Taylor's Theorem, \cite{MR99f:17033}, 2.4.3]
\label{Taylorsthm}  Let
$a(z)$ be a formal distribution.  Then in the region $|z-w|<|w|$,
\begin{equation}
a(z)=\sum_{j=0}^\infty \partial_w^{(j)}a(w)(z-w)^j.
\end{equation}
\end{thm}

\begin{thm}[\cite{MR99f:17033}, Theorem 2.3.2]\label{kac}  Set $\mathbb C[\mathbf x]=\mathbb C[x_n,x^1_n|n\in\mathbb Z]$ and $\mathbb C[\mathbf y]= C[y_m,y_m^1|m\in\mathbb N^*]$.  Let $a(z)$ and $b(z)$
be formal distributions with coefficients in the associative algebra
 $\End(\mathbb C[\mathbf x]\otimes \mathbb C[\mathbf y])$ where we are using the usual normal ordering.   The
following are equivalent
\begin{enumerate}[(i)]
\item
$\displaystyle{[a(z),b(w)]=\sum_{j=0}^{N-1}\partial_w^{(j)}
\delta(z-w)c^j(w)}$, where $c^j(w)\in \End(\mathbb C[\mathbf x]\otimes \mathbb
C[\mathbf y])[\![w,w^{-1}]\!]$.
\item
$\displaystyle{\lfloor
ab\rfloor=\sum_{j=0}^{N-1}\iota_{z,w}\left(\frac{1}{(z-w)^{j+1}}
\right)
c^j(w)}$.
\end{enumerate}\label{Kacsthm}
\end{thm}

In other words the singular part of the {\it operator product
expansion}
$$
\lfloor
ab\rfloor=\sum_{j=0}^{N-1}\iota_{z,w}\left(\frac{1}{(z-w)^{j+1}}
\right)c^j(w)
$$
completely determines the bracket of mutually local formal
distributions $a(z)$ and $b(w)$.   One writes
$$
a(z)b(w)\sim \sum_{j=0}^{N-1}\frac{c^j(w)}{(z-w)^{j+1}}.
$$

\subsection{The elliptic Heisenberg algebra}
The Cartan subalgebra $\mathfrak h$ tensored with $R$ generates a subalgebra of $\hat{{\mathfrak g}}$ which is an extension of an oscillator algebra.    This extension motivates the following definition:  The Lie algebra with generators $b_{m},b_m^1$, $m\in\mathbb Z$, $\mathbf 1_0,\mathbf 1_\pm $, and relations
\begin{align}
[b_{m},b_{n}]&=2n\,\delta_{m+n,0}\mathbf 1_0\label{b1} \\
[b^1_m,b_n^1] &=(n-m) (\delta_{m+n+2,-1}-2b \,\delta_{m+n+2,0} +\delta_{m+n+2,1})\mathbf 1_0 \label{b2}\\
[b^1_m,b_n] &=2n\Big(p_{|m+n+1| }(b) \mathbf 1_+ +q_{|m+n+1|}(b)\mathbf 1_-\Big)\label{b3} \\
[b_{m},\mathbf 1_0]&=[b_{m}^1,\mathbf 1_0]=[b_{m},\mathbf 1_\pm]=[b_{m}^1,\mathbf 1_\pm]= 0.\label{b4}
\end{align}
we will give the appellation the {\it elliptic Heisenberg algebra} and denote it by $\hat{\mathfrak h}$.

This algebra has an involutive anti-automorphism $\sigma$:
$$
\sigma(b_n)=-b_{-n},\quad \sigma(b^1_n)=-b^1_{-n-2},\quad
\sigma(\mathbf 1_0)=\mathbf 1_0,\quad \sigma(\mathbf 1_\pm)=\mathbf 1_\pm.
$$

If we introduce the formal distributions
\begin{equation}
\beta(z):=\sum_{n\in\mathbb Z} b_nz^{-n-1},\quad \beta^1(z):=\sum_{n\in\mathbb Z}b_n^1z^{-n-1}.
\end{equation}
then using calculations done earlier for the elliptic Lie algebra we can see that the relations above can be rewritten in the form
\begin{align*}\label{bosonrelations}
[\beta(z),\beta(w)]&=2\mathbf 1_0\partial_z\delta(z/w)=-2\mathbf 1_0\partial_w\delta(z/w) \\
[\beta^1(z),\beta^1(w)]
&=-2\left(\frac{1}{2}\partial_w(w(w^2 -2bw+1)) \delta(z/w)+w(w^2 -2bw+1)  \partial_w(\delta(z/w)\right) \mathbf 1_0 \\
[\beta^1(z),\beta(w)]
&= -2\left((P(w^{-1},b)+P(w,b))\mathbf 1_++(Q(w^{-1},b)+Q(w,b)-2)\mathbf 1_-\right)w \partial_w\delta(z/w).
\end{align*}

 \color{black}
 Set
\begin{align*}
\hat{\mathfrak h}^\pm:&=\sum_{n\gtrless 0}\left(\mathbb Cb_n+\mathbb Cb_{n-1}^1\right),\quad
\hat{ \mathfrak h}^0:=\mathbb C\mathbf 1_-\oplus \mathbb C\mathbf 1_0\oplus \mathbb C\mathbf 1_+\oplus \mathbb Cb_0\oplus \mathbb Cb^1_{-1}.
\end{align*}
We introduce Borel type subalgebras of $\hat{\mathfrak h}$:
\begin{align*}
\hat{\mathfrak b}= \hat{\mathfrak h}^+ + \hat{\mathfrak h}^0,\qquad \hat{\mathfrak b}^-= \hat{\mathfrak h}^-+ \hat{\mathfrak h}^0.
\end{align*}
The defining relations above ensure that
$\hat{\mathfrak b}$ is a subalgebra. Note that
$\sigma(b^1_{-1})=-b^1_{-1}$, $\sigma(\hat{\mathfrak
h}^0)=\hat{\mathfrak h}^0$,  and $\sigma(\hat{\mathfrak
h}^+)=\hat{\mathfrak h}^-$.  Moreover for $m,n\geq 0$  or for
$m,n<-1$, we get
\begin{align}
[b^1_m,b_n^1] &=(n-m) (\delta_{m+n+2,-1}-2b \,\delta_{m+n+2,0} +\delta_{m+n+2,1})\mathbf 1_0=0 .
\end{align}


\begin{lem}\label{twodimrep}
Let $\mathcal V=\mathbb C\mathbf v_0\oplus \mathbb C\mathbf v_1$ be a two-dimensional representation of $\hat{\mathfrak b}$ where $\hat{\mathfrak h}^+\mathbf v_i=0$ for $i=0,1$.   Suppose  $\lambda,\mu,\nu,\varkappa,\chi_\pm,\chi_0 \in \mathbb C$  are such that
\begin{align*}
b_0\mathbf v_0&=\lambda \mathbf v_0,  &b_0\mathbf v_1&=\lambda \mathbf v_1 \\
b_{-1}^1\mathbf v_0&=\mu \mathbf v_0+\nu \mathbf v_1,  &b_{-1}^1\mathbf v_1&=\varkappa \mathbf v_0+\mu \mathbf v_1\\
\mathbf 1_\pm\mathbf v_i&=\chi _\pm  \mathbf v_i,\quad   &\mathbf 1_0\mathbf v_i&=\chi_0\mathbf v_i,\quad i=0,1.
\end{align*}
Then  $\chi _+=\chi_-=0$.
\end{lem}
\begin{proof} Since $b_m$ acts by scalar multiplication for $m,n\geq 0$, the first defining relation \eqnref{b1} is satisfied for $m,n\geq 0$.
The second relation \eqnref{b2} is also satisfied as the right hand side is zero if $m\geq 0,n\geq 0$.  If $n=0$, then since $b_0$ acts by a scalar, the relation \eqnref{b3} leads to no condition on $\lambda,\mu,\nu,\varkappa,\chi_\pm,\chi_0 \in \mathfrak h_4^0$.    If $m\geq 0$ and $n>0$, the third relation becomes
$$
0=b^1_mb_n\mathbf v-b_nb^1_m\mathbf v=[b^1_m,b_n] \mathbf v=2n\Big(p_{m+n+1 }(b) \mathbf 1_+ +q_{m+n+1}(b)\mathbf 1_-\Big)\mathbf v.
$$
Now
$$
p_{2}(b)=\frac{4}{5}(b-1),\quad q_2(b)=\frac{1}{5},\quad p_{3}(b)=\frac{1}{35} \left(32 b^2-48 b+11\right),\quad q_3(b)=\frac{4}{35} (2 b-1)
$$
which implies that $\chi_+=\chi_-=0$.
\end{proof}
Let $B_{-1}^1$ denote the linear transformation on $\mathcal V$ that agrees with the action of $b_{-1}^1$.
If we define the notion of a $\hat{\mathfrak b}$-submodule as is done in \cite{MR1328538}, Definition 1.2, then $\mathcal V$ above is an irreducible $\hat{\mathfrak b}$-module when $\varkappa \nu\neq 0$ i.e. if $\det B_{-1}^1\neq \mu^2$.  Later we will form  a module for the elliptic affine algebra by creating the induced module for $V$.  The resulting representation can not be irreducible,  if $\mathcal V$ were not irreducible itself (in the sense of Sheinman, [15]).

We would like to give the heuristic construction of formulae for the representations given below.   Since we  do not claim that any of this construction gives us a mathematically rigorous proof that we have obtain a representation, we will still need to check the defining relations are satisfied by the given formulae.

To this end we consider the following expression (which should lie in a Borel subgroup $\hat{B}_-$ of $\hat{\mathfrak b}^-$)
$$
\exp\left(\sum_{m<0} y_m b_m\right)\exp\left(\sum_{m<-1}y^1_mb^1_m\right),
$$
where $y_m, y_n^1$ are coordinate functions. Consider the representation $\mathcal V$ defined above with $\mathbf w\in \mathcal V$.
Since $\exp A\exp B=\exp(\exp(\text{ad} \,A)B)\exp A$ we have for $k> 0$.
\begin{align*}
\exp(-tb_k)&\exp\left(\sum_{m<0} y_m b_m\right)\exp\Big(\sum_{m<-1}y^1_mb^1_m\Big)\mathbf w
         \\
         &=\exp\left(\exp \text{ad}(-tb_k)\sum_{m<0} y_m b_m\right)
         \exp\left(\exp \text{ad}(-tb_k)\sum_{m<-1} y^1_m b^1_m\right) \mathbf w \\
         &=\exp\left( \sum_{m<0} y_m b_m-t\sum_{m<0} y_m [b_k,b_m]\right)
         \exp\left(\sum_{m<-1} y^1_m b^1_m-t\sum_{m<-1} y^1_m [b_k,b^1_m]\right) \mathbf w +\mathcal O(t^2)\\
         &=\exp\left( \sum_{m<0} y_m b_m-t\sum_{m<0} y_m 2m\delta_{k+m,0}\mathbf 1_0\right)\times \\
         &\exp\left(\sum_{m<-1} y^1_m b^1_m+t\sum_{m<-1} y^1_m 2k\Big(p_{|m+k+1| }(b) \mathbf 1_+ +q_{|m+k+1|}(b)\mathbf 1_-\Big) \right) \mathbf w +\mathcal O(t^2).
\end{align*}
Under the assumption $[A,A,B]=0$ we can use Campbel-Baker-Hausdorff formula:
$$
\log (\exp (A+tB)\exp(  -tB)=A-\frac{t}{2}[A,B]+O(t^2).
$$
This leads to
\begin{align*}
\exp(-tb_k)&\exp\left(\sum_{m<0} y_m b_m\right)\exp\Big(\sum_{m<-1}y^1_mb^1_m\Big)\mathbf w
         \\
         &=\exp\left( \sum_{m<0} y_m b_m\right)\exp\left(\sum_{m<-1} y^1_m b^1_m\right) \times \\
         &\exp\left(-t\sum_{m<0} y_{m} 2m\delta_{k+m,0}\mathbf 1_0\right) \exp\left(t\sum_{m<-1} y^1_m 2k\Big(p_{|m+k+1| }(b) \mathbf 1_+ +q_{|m+k+1|}(b)\mathbf 1_-\Big) \right) \mathbf w +\mathcal O(t^2) \\
         &=\exp\left( \sum_{m<0} y_m b_m\right)\exp\left(\sum_{m<-1} y^1_m b^1_m\right) \exp\left(t\,2k y_{-k} \mathbf 1_0\right) \mathbf w +\mathcal O(t^2).
\end{align*}
Thus we set
$$
\rho(b_k)=-2k\chi_0\,y_{-k}.
$$

After conjugating the linear maps $\rho$ (defined at the moment only on the $b_k$ with $k>0$) with the
anti-automomorphism $\sigma$ together with the ``anti-Fourier
transform''
$$
\Phi: y_k\mapsto -\partial_{y_k},\quad \partial_{y_k}\mapsto -y_k,
$$
we get a new linear map which we also denote by $\rho$:
$$
\rho(b_{-k})=-2k \chi_0\,\partial_{y_{-k}}, \,\, k>0.
$$
 The formula for $\rho$ on the other basis elements of $\hat{\mathfrak h}$ are obtained in a similar fashion and are given in the following

\begin{prop}\label{heisenbergprop}
Let $\mathcal M=\mathbb C[y_{-n},y_{-m}^1|\,m,n\in\mathbb N^*]\otimes \mathcal V$.  Then for $k>0$
\begin{align*}
\rho(b_{k})&=-y_{-k},\quad \rho(b_{-k-1}^1)=y_{-k-1}^1    \\
\rho(b_{-k})&=-2k \chi_0\,\partial_{y_{-k}, }
 \\
\rho(b_{k}^1)&=-\left((2k+1) \partial_{y_{-k-1}^1}  -4(k+1) b\,
\partial_{y_{-k-2}^1}   +(2k+3)
\partial_{y_{-k-3}^1}\right) \chi_0     \\
\rho(b_{0}^1)&=-\left(\partial_{y_{-1}^1}  -4b\,
\partial_{y_{-2}^1}   +3
\partial_{y_{-3}^1}\right) \chi_0     \\
\rho(b_{-1}^1)&=y_{-1}^1-\partial_{y_{-2}^1}\chi_0+B_{-1}^1  \\
\rho(b_0)&=\lambda\quad   \rho(\mathbf 1_0)=\chi_0,\quad \rho(\mathbf 1_\pm)=0
\end{align*}
defines a representation of $\hat{\mathfrak h}$ on $\mathcal M$.
\end{prop}

\begin{proof}
The defining relations on the generators of $\hat{\mathfrak h}$ can easily be checked. For example,
for $m\geq 0$ and $n>0$ we have
\begin{align*}
[&\rho(b_{m}),\rho(b_{-n})]\\
&=[-y_{-m}, -2n\chi_0\partial_{y_{-n}}]   \\
&=-2n\,\delta_{m-n,0}\rho(\mathbf 1_0),
\end{align*}
or

for $n>1$ and $m\geq 0$,
\begin{align*}
[&\rho(b_{-n}^1 ),\rho(b_{m}^1 )]\\
&=\Big[y_{-n}^1,
  -(2m+1)\partial_{y_{-m-1}^1} + 4(m+1) b\,\partial_{y_{-2-m}^1} -(2m+3)\partial_{y_{-m-3}^1}   )\chi_0 \Big]  \\
 &=-\left(-(2m+1)\delta_{-n+m,-1}  +4(m+1) b\,\delta_{-n+m,-2}  -(2m+3)\delta_{-n+m,-3}  \right)\chi_0  \\
&=(m+n)\left(\delta_{-n+m,-1} -2b\,\delta_{-n+m,-2}
+\delta_{-n+m,-3}  \right)\chi_0 .
\end{align*}

All other relations are straightforward.
\end{proof}

\begin{subsection}{The Imaginary Borel subgroup}\label{borel}


Let $e,h,f$ be the standard basis of $sl(2)$. Then $e_n$, $e_n^1$,
$f_n$, $f_n^1$, $h_n$, $h_n^1$, $w_0$, $w_{\pm}$ is a basis of
$\hat{\mathfrak g}$, where $x_n=x\otimes t^n$, $x_n^1=x\otimes
ut^n$.  Denote by $\mathfrak N_{+}$ (respectively $\mathfrak
N_{-}$) a subalgebra spanned by $e_n$, $e_n^1$, $n\in \mathbb Z$,
$h_m$, $h_m^1$, $m>0$ (respectively $f_n$, $f_n^1$, $n\in \mathbb
Z$, $h_m$, $h_m^1$, $m<0$) and set $\mathfrak H=\mathbb C h\oplus
\mathbb C w_0\oplus \mathbb C w_+ \oplus \mathbb C w_-$. Then
$$\hat{\frak g}=\mathfrak N_{-}\oplus  \mathfrak H \oplus \mathfrak N_{+}.$$

We call the subalgebra $\mathfrak B=\mathfrak N_+\oplus H$ the
Imaginary Borel subalgebra of $\hat{\mathfrak g}$. We will make
use of this Borel subalgebra to construct Wakimoto type
representations of $\hat{\mathfrak g}$.
%


\end{subsection}

\section{Two realizations of the elliptic affine algebra $\hat{{\mathfrak g}}$}
Our main result is the following
\begin{thm} \label{mainthm}  Let $r\in\{0,1\}$, with the corresponding normal ordering defined above.
Assume that $\chi_+=\chi_-=0$, $\chi_0\in\mathbb C$ and $\mathcal
V$ as in \propref{heisenbergprop}.   Then the following defines a
representation of the elliptic affine salgebra $\hat{\mathfrak g}$ on
$\mathbb C[\mathbf x]\otimes \mathbb C[\mathbf y]\otimes \mathcal
V$:
\begin{align*}
\theta(\omega_\pm)&=0, \qquad
\theta(\omega_0)=\chi_0=\chi_0-4\delta_{r,0}   \\
\theta(f(z))&=-\alpha, \qquad
\theta(f^1(z))= -\alpha^1,   \\ \\
\theta(h(z))
&=2\left(:\alpha\alpha^*:+:\alpha^1\alpha^{1*}: \right)
     +\beta   \\  \\
\theta(h^1(z))
&=2\left(:\alpha^1\alpha^*:
   +z(1-2bz+z^2):\alpha\alpha^{1*}: \right)
     +\beta^1  \\  \\
\theta(e^1(z))
&=:\alpha^1(z)\alpha^*(z)^2:
        +z\left( z^{2} -2b z+ 1\right)\left(  :\alpha^1(z)\alpha^{1*}(z)^2:   +2  :\alpha(z)\alpha^*(z)\alpha^{1*}(z):\right) \\
&\quad +\beta^1(z)\alpha^*(z)+z\left( z^{2} -2b z + 1\right)\beta(z)\alpha^{1*}(z)\\
&\quad+\chi_0z(1-2bz+z^2)\partial_z a^{1*}(z)+\frac{1}{2}\chi_0\partial_z[z(1-2bz+z^2)] a^{1*}(z)\\   \\
\theta(e(z))
&=:\alpha(\alpha^*)^2: +z(1-2bz+z^2):\alpha(\alpha^{1*})^2: +2 :\alpha^1\alpha^*\alpha^{1*}: +\beta\alpha^*  \\
&\quad +\beta^1\alpha^{1*} +\chi_0\partial\alpha^*.
\end{align*}
\end{thm}


In the formulas above we omitted the variables in the fields
whenever it does not create any confusion.   In the next section we explain in what sense $r=0$ and $r=1$ give two different realizations of the elliptic affine Lie algebra $\hat{\mathfrak g}$.

Before we go through the proof it will be fruitful to introduce
V.Kac's $\lambda$-notation (see \cite{MR99f:17033} section 2.2 and
\cite{MR1873994} for some of its properties) used in operator
product expansions.  If $a(z)$ and $b(w)$ are formal
distributions, then
$$
[a(z),b(w)]=\sum_{j=0}^\infty \frac{(a_{(j)}b)(w)}{(z-w)^{j+1}}
$$
is transformed under the {\it formal Fourier transform}
$$
F^{\lambda}_{z,w}a(z,w)=\text{Res}_ze^{\lambda(z-w)}a(z,w),
$$
 into the sum
\begin{equation*}
[a_\lambda b]=\sum_{j=0}^\infty \frac{\lambda^j}{j!}a_{(j)}b.
\end{equation*}


\begin{proof}
We have need to check that the following table is preserved under  $\theta$.
\begin{table}[htdp]
\caption{}
\begin{center}
\begin{tabular}{c|cccccc}
$[\cdot_\lambda \cdot]$ & $f(w)$ & $f^1(w)$ & $h(w)$ & $h^1(w)$ & $e(w)$ & $e^1(w)$ \\
\hline
$f(z)$ & $0$ &  $0$ & $*$ &  $*$ &$*$  &$*$    \\
$f^1(z)$ &   &  $0$  & $*$  & $*$  & $*$  &$*$  \\
$h(z)$ &  &   &  $*$  & $*$  & $*$  & $*$ \\
$h^1(z)$ &&   &     & $*$  &    $*$ & $*$ \\
$e(z)$  &&   &     &  &    $0$ &  $0$  \\
$e^1(z)$ & &   &     &  &   &  $0$ \\
\end{tabular}
\end{center}
\label{default}
\end{table}%

Here $*$ represent nonzero formal distributions that are obtained
from the the defining relations. The proof is carried out using
Wick's Theorem together with Taylor's Theorem.

\begin{align*}
[\theta(f)_\lambda\theta(f)]&=0,\quad [\theta(f)_\lambda\theta(f^1)]= 0,\\ \\
[\theta(f)_\lambda\theta(h)]&=-\Big[\alpha_\lambda\Big(2\left(\alpha\alpha^*+\alpha^1\alpha^{1*}
\right)   +\beta \Big)\Big]=
- 2\alpha =  - 2\theta(f) , \\  \\
[\theta(f)_\lambda\theta(h^1)]&=-[\alpha_\lambda\left(2\left(\alpha^1\alpha^*
   +P^2\alpha\alpha^{1*} \right)
     +\beta^1\right)  ]= -2\alpha^{1}=   2\theta(f^1), \\  \\
[\theta(f)_\lambda\theta(e)]&=
-[\alpha_\lambda\Big(:\alpha(\alpha^*)^2:+\beta\alpha^* +
\chi_+\partial ( P  \alpha^{1*})
        +P^2:\alpha(\alpha^{1*})^2:\\
&\hskip 100pt +2 :\alpha^1\alpha^*\alpha^{1*}: +\beta^1\alpha^{1*} +\chi_0\partial\alpha^*\Big)] \\
&=-2\left(:\alpha\alpha^*:    - :\alpha^1\alpha^{1*}:
\right)-\beta-\chi_0\partial
= -\theta(h)-\chi_0\lambda   \\ \\
[\theta(f)_\lambda\theta(e^1)]&=- \Big[\alpha_\lambda\Big(:\alpha^1(\alpha^*)^2: +\beta^1\alpha^*+\chi_0P(\partial P)\alpha^{1*}  + \chi_+ P\partial\alpha^*  \\
&\qquad
        +P^2\left(  :\alpha^1(\alpha^{1*})^2:   +2  :\alpha\alpha^*\alpha^{1*}:
        +\beta\alpha^{1*}+\chi_0\partial\alpha^{1*}\right) \Big)\Big]\\
&=-2\left(:\alpha^1\alpha^*:    + P^2:\alpha\alpha^{1*}: \right)-\beta^1 -\chi_+ P  \partial =-\theta(h^1)-\chi_+P \lambda
\\\\ \\  \\
 [\theta(f^1)_\lambda\theta(f^1)]&=0  \\ \\
[\theta(f^1)_\lambda\theta(h)]&=-[\alpha^1_\lambda\left(2\left(:\alpha\alpha^*:+:\alpha^1\alpha^{1*}: \right)   +\beta\right) ]=  -2\alpha^1 =   2\theta(f^1),  \\  \\
[\theta(f^1)_\lambda\theta(h^1)]&=-[\alpha^1_\lambda\left(2\left(:\alpha^1\alpha^*:
   +P^2:\alpha\alpha^{1*}: \right)
     +\beta^1\right)] =-2P^2\alpha^{1} =P^2\theta(f^{1}) , \\  \\
[\theta(f^1)_\lambda\theta(e)]&=- [\alpha^1_\lambda\Big(:\alpha(\alpha^*)^2:+\beta\alpha^* + \partial (P\alpha^{1*})   +P^2:\alpha(\alpha^{1*})^2: +2 :\alpha^1\alpha^*\alpha^{1*}: +\beta^1\alpha^{1*} +\chi_0\partial\alpha^*\Big)] \\
&=- \Big(  \chi_+(\partial P+P\lambda ) +2P^2:\alpha\alpha^{1*}: +2 :\alpha^1\alpha^*: +\beta^1  \Big)
=-\theta(h^1)-\chi_+(\partial P+P\partial) \\ \\
[\theta(f^1)_\lambda\theta(e^1)]&=-[ \alpha^1_\lambda\Big(:\alpha^1(\alpha^*)^2: +\beta^1\alpha^*+\chi_0P(\partial P)\alpha^{1*}  + \chi_+P \partial\alpha^*  \\
&\qquad
        +P^2\left(  :\alpha^1(\alpha^{1*})^2:   +2  :\alpha\alpha^*\alpha^{1*}:
        +\beta \alpha^{1*}+\chi_0\partial\alpha^{1*}\right) \Big)] \\
&=-\Big(\chi_0P(\partial P)\alpha^{1*}
        +P^2\left(  2\left(:\alpha^1\alpha^{1*}:   + :\alpha\alpha^* :\right)
        +\beta+\chi_0\lambda  \right) \Big)  \\
&=-\left(P^2\theta(h)+\chi_0P\partial P+P^2\chi_0\lambda\right)
\end{align*}

\begin{align*}
[\theta(h)_\lambda\theta(h)]&=\left(2\left(:\alpha\alpha^*:+:\alpha^1\alpha^{1*}: \right)   +\beta\right)  _\lambda \left(2\left(:\alpha\alpha^*:+:\alpha^1\alpha^{1*}: \right)   +\beta(w)\right)]  \\
&= 4\Big(-:\alpha\alpha^*: +:\alpha^*\alpha: -:\alpha^1\alpha^{1*}: +:\alpha^{1*}\alpha^1:  \Big)
-8\delta_{r,0}\lambda +[\beta_\lambda\beta]   \\
&=-2(4\delta_{r,0}+\chi_0)\lambda
 \end{align*}
 Thus
 \begin{align*}
\left[\theta(h(z)),\theta(h(w))\right]&=
 -2(4\delta_{r,0}+\chi_0)\partial_{w}\delta(z/w).
  \end{align*}
Next we calculate
\begin{align*}
[\theta(h)_\lambda\theta(h^1)]&=\left[\left(2\left(:\alpha\alpha^*:+:\alpha^1\alpha^{1*}: \right)   +\beta \right) _\lambda \left(2\left(:\alpha^1\alpha^*:
   +P^2:\alpha\alpha^{1*}: \right)+\beta^1  \right) \right]\\
&=4\Big(\left(:\alpha^* \alpha^1 : -:\alpha^1 \alpha^{*} :\right)
+ P^2 \left(-:\alpha\alpha^{1*} : +:\alpha^{1*}\alpha :  \right)
\Big)  +[\beta_\lambda\beta^1].
\end{align*}
Since $[a_n,a_m^{1*}]=[a^1_n,a_m^{*}]=0$, we have
 \begin{align*}
\left[\theta(h(z)),\theta(h^1(w))\right]&=[\beta(z),\beta^1(w)]=-2
\chi_+\left(1-2bw+w^2\right)^{1/2}\partial_w\delta(z/w).
 \end{align*}
 We continue with
\begin{align*}
[\theta(h^1)_\lambda\theta(h^1)]&=\left[2\left(:\alpha^1\alpha^*:
   +P^2:\alpha\alpha^{1*}: \right)
     +\beta^1  \right)_\lambda \left(2\left(:\alpha^1 \alpha^* :
   +P^2:\alpha\alpha^{1*}: \right)
     +\beta^1 )  \right] \\
&= 4P^2\left(-:\alpha\alpha^*:+:\alpha^{1*}\alpha^1:\right) +4P^2\left(-:\alpha^1\alpha^{1*} : +:\alpha^*\alpha: \right)
  \\
&\quad  +8\delta_{r,0}P^2\lambda   +8\delta_{r,0}P(\partial P)    + [\beta^1_\lambda \beta^1]   \\
&=8\delta_{r,0}P^2\lambda   +8\delta_{r,0}P(\partial P)    -2\chi_0(P^2\lambda+P(\partial P)).
 \end{align*}
 Note that $:a(z)b(z):$ and $:b(z)a(z):$ are usually not equal, but
 $:\alpha^1(w)\alpha^{*1}(w):=:\alpha^{1*}(w)\alpha^1(w):$ and $:\alpha(w)\alpha^*(w):=:\alpha^*(w)\alpha (w):$.   Thus we have
\begin{align*}
\left[\theta(h^1(z)),\theta(h^1(w))\right]&=(4\delta_{r,0}+\chi_0)\left(-2\left(w^2-2bw+1\right)\partial_w\delta(z/w)
+2  (b-w)\delta(z/w)\right).
\end{align*}

Next we calculate the $h$'s paired with the $e$'s:

\begin{align*}
[\theta(h)_\lambda \theta(e^1)]&=\Big[\Big(2\left(:\alpha\alpha^*:+:\alpha^1\alpha^{1*}: \right)
     +\beta \Big)_\lambda \\
&\quad     \Big(:\alpha^1(\alpha^*)^2: +\beta^1\alpha^*+\chi_0P(\partial P)\alpha^{1*}
 +\chi_+ P \partial\alpha^*  \\
&\quad
        +P^2\left(  :\alpha^1(\alpha^{1*})^2:   +2  :\alpha\alpha^*\alpha^{1*}:
        +\beta\alpha^{1*}+\chi_0\partial\alpha^{1*}\right)\Big)\Big] \\  \\
&=2:\alpha^1(\alpha^*)^2:
            +2\beta^1\alpha^*
        +2\chi_0P(\partial P)\alpha^{1*}  +2\chi_+ P \alpha^*\lambda +2\chi_+ \partial P  \\
&\quad
+4P^2\Big(:\alpha^*\alpha\alpha^{1*}:
    - 2\delta_{r,0}\alpha^{1*} \lambda \Big)
        -2P^2\Big(  :\alpha^1(\alpha^{1*})^2:  +\beta\alpha^{1*} + \chi_0\lambda\Big) \\
&\quad+  \alpha^* [\beta_\lambda\beta^1]+P^2\alpha^{1*} [\beta_\lambda \beta]
\end{align*}
From this we can conclude
\begin{align*}
\left[\theta(h(z)),\theta(e^1(w))\right]
&=2\left(:\alpha^1(w)\alpha^*(w)^2:+\beta^1(w)\alpha^*(w)
        +\chi_0(w-b)\alpha^{1*}(w)\right)\delta(z/w)  \\
&\quad +2\left(\chi_+ (1-2bw+w^2)^{1/2} \partial_w\alpha^*(w) \right)\delta(z/w) \\
&\quad
        +2\left( w^{2} -2b w  + 1\right)\Big(2:\alpha(w)\alpha^*(w)\alpha^{1*}(w):
        +\chi_0\partial_w\alpha^{1*}(w)\Big)\delta(z/w) \\
&\quad
        +2\left( w^{2} -2b w  + 1\right)\Big( :\alpha^1(w)\alpha^{1*}(w)\alpha^{1*}(w)
                +\beta(w)\alpha^{1*}(w)\Big)\delta(z/w) \\  \\
&\quad -2\left( w^{2} -2b w  + 1\right)(4\delta_{r,0}+ \chi_0 -\chi_0) \alpha^{1*}(w) \partial_w \delta(z/w)\\  \\
&=2\theta(e^1(w))\delta(z/w)  -2\left( w^{2} -2b w  + 1\right)(4\delta_{r,0}+ \chi_0 -\chi_0) \alpha^{1*}(w) \partial_w \delta(z/w).
\end{align*}

Next we must calculate

\begin{align*}
[\theta(h)_\lambda\theta(e)]
&=2[:\alpha\alpha^*:
_\lambda   \Big(:\alpha(\alpha^*)^2:+\beta\alpha^* + \chi_+\partial P\alpha^{1*}
   +\chi_+P  \partial\alpha^{1*}\\
&\quad +P^2:\alpha(\alpha^{1*})^2: +2 :\alpha^1\alpha^*\alpha^{1*}: +\beta^1\alpha^{1*} +\chi_0\partial\alpha^*\Big) ]\\  \\
&\quad+2[:\alpha^1\alpha^{1*}: _\lambda \Big(:\alpha(\alpha^*)^2:+\beta\alpha^* + \chi_+(\partial P)\alpha^{1*}  +\chi_+ P  \partial\alpha^{1*}\\
&\quad +P^2:\alpha(\alpha^{1*})^2: +2 :\alpha^1\alpha^*\alpha^{1*}: +\beta^1\alpha^{1*} +\chi_0\partial \alpha^*
\Big)] \\  \\
&\quad
     +[\beta _\lambda \Big(:\alpha(\alpha^*)^2:+\beta\alpha^* +\chi_+(\partial P)\alpha^{1*}
  +\chi_+P \partial\alpha^{1*}\\
&\quad +P^2:\alpha(\alpha^{1*})^2: +2 :\alpha^1\alpha^*\alpha^{1*}: +\beta^1\alpha^{1*} +\chi_0\partial_w\alpha^*
\Big) ]\\  \\
&=2 \Big(:\alpha(\alpha^*)^2: +\beta\alpha^*
+ \chi_+\partial (P\alpha^{1*})  \\
&\quad +P^2:\alpha(\alpha^{1*})^2:+2 :\alpha^1\alpha^*\alpha^{1*}:  +\beta^1 \alpha^{1*} +\chi_0\partial\alpha^*\Big) \\
&\quad 2\chi_0\alpha^*\lambda-8\delta_{r,0}\alpha^*\lambda
-2\chi_0\alpha^*\lambda.\end{align*}

As a consequence we get
\begin{align*}
\left[\theta(h(z)),\theta(e(w))\right]&
=2 \theta(e(w))\delta(z/w)+2(\chi_0-\chi_0-4\delta_{r,0}) \alpha^*(w) \partial_w\delta(z/w).
\end{align*}

Next we calculate
\begin{align*}
[\theta(h^1)_\lambda\theta(e^1)]&=\Big[\Big(2\left(:\alpha^1\alpha^*:
   +P^2:\alpha\alpha^{1*}: \right)
     +\beta^1  \Big)_ \lambda   \Big(:\alpha^1(\alpha^*)^2: +\beta^1\alpha^*+\chi_0P(\partial P)\alpha^{1*}
     +\chi_+ P \partial\alpha^*  \\
&\quad
        +P^2\left(  :\alpha^1(\alpha^{1*})^2:   +2  :\alpha\alpha^*\alpha^{1*}:
        +\beta\alpha^{1*}+\chi_0\partial\alpha^{1*}\right)\Big)\Big] \\  \\
&=      2  P^2:\alpha(\alpha^*)^2:
                +2P^2\beta\alpha^{*}+2\chi_0P^2\partial\alpha^{*} +2P^4:\alpha(\alpha^{1*})^2: \\
&+  4P^2:\alpha^1\alpha^{1*}\alpha^*:    +2P^2\beta^1\alpha^{1*} +2\chi_+ P^3 \partial\alpha^{1*}
         +2\chi_+ P^2(\partial P)\alpha^{1*} \\  \\
& +2\chi_0\alpha^* (P(\partial P)  + P^2\lambda)   -8\delta_{r,0}\alpha^*(P(\partial P)+P^2\lambda) \\
&-2\chi_0 \alpha^* (P(\partial P)+P^2\lambda) \\ \\
&=2P^2\theta(e).
\end{align*}

The final calculation for the Cartan generators is
\begin{align*}
[\theta(h^1)_\lambda \theta(e)]
&=\Big[\Big(2\left(:\alpha^1\alpha^*:+P^2:\alpha\alpha^{1*}: \right)
     +\beta^1 \Big)_\lambda   \Big(:\alpha (\alpha^*)^2:+\beta \alpha^*  +  \chi_+(\partial P)\alpha^{1*}
  +\chi_+ P  \partial \alpha^{1*}\\
&\hskip 150pt +P^2:\alpha (\alpha^{1*})^2: +2 :\alpha^1 \alpha^* \alpha^{1*} : +\beta^1 \alpha^{1*} +\chi_0\partial \alpha^*
\Big) \Big]\\  \\
&= -2:\alpha^1 (\alpha^*)^2  + 2 \chi_+(\partial P)\alpha^{*}
  +2\chi_+ P  \partial \alpha^{*}+2\chi_+ P  \alpha^{*}\lambda \\
&\qquad  -2P^2:\alpha^1 (\alpha^{1*})^2: +4P^2:\alpha \alpha^*\alpha^{1*}: -4\delta_{r,0} P^2 \alpha^{1*}\lambda
+ 4 :\alpha^1 (\alpha^*)^2 : +\beta^1 \alpha^{*}
 \\  \\
&\quad+4P^2:\alpha \alpha^*\alpha^{*1}: +2P^2 \beta \alpha^{*1} +4P^2 :\alpha^1 ( \alpha^{1*})^2 : -4P^2 :\alpha\alpha^* \alpha^{1*} : \\
&\quad
-4 \delta_{r,0}P^2 \alpha^{1*} \lambda -8 \delta_{r,0}P(\partial P) \alpha^{1*}  +2\chi_0P^2\alpha^{*1}\lambda
+2\chi_0P^2\partial \alpha^{1*}  +4\chi_0P(\partial P) \alpha^{1*}  \\  \\
&\quad   -2\chi_+(P\lambda+ \partial P)\alpha^*   -2\chi_0 \alpha^{1*}(P(\partial P)+ P^2 \lambda)      \\  \\
&=2 :\alpha^1 (\alpha^*)^2 +  2\chi_+ P  \partial \alpha^{*}   +2P^2:\alpha^1 (\alpha^{1*})^2: +
 +\beta^1 \alpha^{*} +4P^2:\alpha \alpha^*\alpha^{*1}:
 \\
&\quad+2P^2 \beta \alpha^{*1}  +2\chi_0P^2\partial \alpha^{1*}  +2\chi_0P(\partial P) \alpha^{1*}   \\ \\
&\quad +2(\chi_0 -4\delta_{r,0}-\chi_0)\alpha^{*1}(P(\partial P)+
P^2 \lambda)=2\theta(e^1).
\end{align*}

Similar calculations show that $$[\theta(e)_\lambda
\theta(e)]=[\theta(e)_\lambda\theta(e^1)]=[\theta(e^1)_\lambda
\theta(e^1)]=0.$$ We leave the details to the reader. This
completes the proof of the theorem.
\end{proof}

\section{Jakobsen-Kac Realizations}\label{JKsection}

Now we relate the representations constructed in
Theorem~\ref{mainthm} with the representations constructed by
H.Jakobsen and V.Kac \cite{JK}.

 If $R$ is an associative (but not necessarily commutative) algebra over $\mathbb C$, then a linear map $\phi:R\to \mathbb C$ is called a {\it trace} on $R$ if
 $$
 \phi(ab)=\phi(ba)\quad\forall a,b\in R.
 $$
 If $\phi$ is a trace on $R$, then the space of matrices
 $$
 \mathfrak{sl}_2(R):=\left\{\begin{pmatrix} a  &  b \\ c & d\end{pmatrix} \,\Big| \,a,b,c,d\in R,\text{ and }\phi(a+d)=0\right\}
 $$
 is a Lie algebra under the usual commutator.

Note that $
 \mathfrak{sl}_2(R)$ always contains a subalgebra isomorphic to $\mathfrak{sl}(2,
 R)=
 \mathfrak{sl}_2\otimes R$.

 The subspace
 $$
 \mathfrak s:=\left\{\begin{pmatrix} a  & b \\  0 & d\end{pmatrix} \,\Big| \,a,b,d\in R,\text{ and }\phi(a+d)=0\right\}
 $$
 is also clearly a subalgebra of $\mathfrak{sl}_2(R)$.  Let $\mathbb C_\phi=\mathbb C\mathbf v_\phi$ be the one
 dimensional $\mathfrak s$-module defined by
 $$
 \begin{pmatrix} a  & b \\  0 & d\end{pmatrix}\mathbf v_\phi:=\phi(a)\mathbf v_\phi.
 $$

  Thus one can define an induced module
 $$
 M(\phi):=U( \mathfrak{sl}_2(R))\otimes_{U(\mathfrak s)}C_\phi.
 $$

 \color{black}  If $R$ is a commutative ring with basis $\{a_\beta\}_{\beta\in I}$ with structure constants $c_{\alpha\beta}^\gamma$ so that
 $$
 a_\alpha a_\beta=c_{\alpha\beta}^\gamma a_\gamma
 $$
  then elements of the form
  $$
  (a_{\alpha_1} f)\cdots   (a_{\alpha_n} f).v_\phi
  $$
constitute a basis of $M(\phi)$.  Here
$\displaystyle{af:=\begin{pmatrix} 0 & 0 \\ a & 0 \end{pmatrix}}$.

 In \cite{JK} Jakobsen and Kac realized  $M(\phi)$ of $\mathfrak{sl}_2(R)$ on the space of
 polynomials $\mathbb C[x_\beta|\beta \in I]$ as follows. Let
 $\tilde{\rho}$ be the corresponding representation. We have

 \begin{align*}
 \rho(a_{\alpha_0}f)&=x_{\alpha_0},  \\
\rho(a_{\alpha_0}h)&= -2  c_{\alpha_0\alpha}^\gamma x_\gamma \frac{\partial}{\partial x_\alpha}+\phi(a_{\alpha_0}) \\
\rho(a_{\alpha_0}e)&=- c_{\alpha_0 \alpha}^\gamma c_{\gamma\beta}^\delta x_\delta \frac{\partial}{\partial x_\alpha}\frac{\partial}{\partial x_\beta}+\phi(a_\gamma)c_{\alpha_0 \alpha}^\gamma \frac{\partial}{\partial x_\alpha}
 \end{align*}

 Now consider $R=\mathbb C[t,t^{-1},u|u^2=t^3-2bt^2+t]$ with  a basis
$\{a_n=t^n, a^1_n:=t^nu\,|\,n\in\mathbb Z\}$. Then
 $$
 a_ma_n=a_{m+n},\quad a_ma^1_n=a_{m+n}^1,\quad a_m^1a_n^1
 =a_{m+n+3}-2ba_{m+n+2}+a_{m+n+1}.
 $$
  The representation $\tilde{\rho}$ of $\mathfrak{sl}_2(R)$ on $\mathbb C[x_n,x_m^1|\,m,n\in\mathbb Z]$ induces a
  representation of $\mathfrak{sl}(2, R)$  given
   by
 \begin{align*}
 \rho(t^m\,f)&=x_m,\quad  \rho(t^mu\,f)=x_m^1 \\
\rho(t^m\, h)&=-2 \sum_p \left(x_{m+p} \partial_{x_p}+x_{m+p}^1  \partial_{x_p^1}\right)+\phi(t^m)  \\
\rho(t^mu\,h)&=-2 \sum_p \left(x_{m+p}^1
\partial_{x_p}+\left(x_{m+p+2}-2x_{m+p+1}+x_{m+p}\right)
\partial_{x_p^1}\right)
+\phi(t^mu)  \\
\rho(t^m\, e)&=-\sum_{n,q}\left(x_{m+n+q}\partial_{x_n}\partial_{x_q}+x_{m+n+q}^1\partial_{x_n}\partial_{x_q^1}+
x_{m+n+q}^1\partial_{x_n^1}\partial_{x_q} \right) \\
&\quad-\sum_{n,q}\left(x_{m+n+q+2}\partial_{x_n^1}\partial_{x_q^1}-2bx_{m+n+q+1}\partial_{x_n^1}\partial_{x_q^1}+x_{m+n+q}
\partial_{x_n^1}\partial_{x_q^1} \right) \\
&\quad  +\sum_p \left(\phi(t ^{m+p} )\partial_{x_p}+ \phi(t ^{m+p}u)\partial_{x_p^1}\right) \\
\rho(t^mu\, e)&= -\sum_{n,q} x_{m+n+q}^1\partial_{x_n}\partial_{x_q } \\
&\quad-\sum_{n,q}\left(x_{m+n+q+2} -2bx_{m+n+q+1}+x_{m+n+q}
  \right)\partial_{x_n}\partial_{x_q^1} \\
&\quad-\sum_{n,q}\left(x_{m+n+q+2}-2bx_{m+n+q+1}+x_{m+n+q} \right)\partial_{x_n^1}\partial_{x_q} \\
&\quad-\sum_{n,q}\left(x_{m+n+q+2}^1-2bx_{m+n+q+1}^1+x_{m+n+q}^1 \right)\partial_{x_n^1}\partial_{x_q^1}\\
&\quad +\sum_p \left(\phi(t^{m+p} u)\partial_{x_p}+ \phi(t
^{m+p}(t^2-2bt-1) )\partial_{x_p^1}\right).
 \end{align*}

One can see that, up to a change in sign, the Jakobsen and Kac's
representation is  a quotient of the representation that we have
constructed in \thmref{mainthm} for the universal central
extension of $\mathfrak{sl}(2,R)$ when $r=1$ and $\phi=0$.  For example if one looks at the Fourier modes of
$$
:\alpha(\alpha^*)^2: +z(1-2bz+z^2):\alpha(\alpha^{1*})^2: +2 :\alpha^1\alpha^*\alpha^{1*}: +\beta\alpha^*
$$
in $\theta(e(z))$, one gets $\rho(t^m\, e)$ up to a difference in sign.

When $r=0$ the Wakimoto type module constructed in
\thmref{mainthm} appears to be related to the Verma module for
$\tilde{\mathfrak g}$ with highest weight subspace $\mathcal V$ given in \lemref{twodimrep}.
We conjecture that generically the Wakimoto type module
constructed in \thmref{mainthm} ($r=0$) and the corresponding
Verma module are isomorphic.



\section{Acknowledgements}
The second author would like to thank the College of Charleston for his fall 2007 sabbatical and also the FAPESP and the University of S\~ao Paulo for the support for his travel  in Brazil during part of this fall  semester.
The third author was partially supported by Fapesp (processo 2005/60337-2) and CNPq (processo 301743/2007-0).


%

\def\cprime{$'$} \def\cprime{$'$}
\providecommand{\bysame}{\leavevmode\hbox to3em{\hrulefill}\thinspace}
\providecommand{\MR}{\relax\ifhmode\unskip\space\fi MR }
\providecommand{\MRhref}[2]{%
  \href{http://www.ams.org/mathscinet-getitem?mr=#1}{#2}
}
\providecommand{\href}[2]{#2}

 \end{document}